\documentclass[12pt]{article}
\date{}
\usepackage{graphicx}
\begin{document}
\title{The well-ordered (F) spaces are $D$-spaces}
\author{{X\small{u} \large{Y}\small{uming}} \thanks{\footnotesize This work is supported by NSF of Shandong Province(No.ZR2010AM019)} \ \ \
\\
{\small $School \ of \ Mathematics, \ Shandong \ University,\ Jinan,\ China$}\\
 }
\maketitle

\hspace{5mm} \begin{minipage}{12.5cm}
{\small \bf \hspace{0cm}Abstract:} {\small We studied the relationships between Collins-Roscoe mechanism and
$D$-spaces, proved that well-ordered $(F)$ spaces are $D$-spaces. This improved the previous results written by D.Soukup and Y.Xu}

\par {\small \bf  \hspace{0cm}Keywords:}
{\small Collins-Roscoe mechanism; well-ordered $(F)$; $D$-space}

\par {\small \bf  \hspace{0cm}AMS classification:}
{\small Primary 54D15; 54D20. Secondary 54E20}

\end{minipage}\\ \\

\begin{center}{\bf \Large 1. Introduction}\end{center}

\vspace{5mm}
In [1], Collins and Roscoe first introduced their well-known structuring mechanism, which was abstracted from the standard proof of the fact that a separable metric space is second countable and has been proved  to be a flexible tool for studying generalized metric spaces.
Let $X$ be a space and for each $x \in X$, $\mathcal{W}(x)= \{W(n, x): n \in \omega \}$ a family of subsets of $X$ containing
$x$. We say that $X$ satisfies condition $(G)$ if,

\hspace{2mm} given open set $U$ containing $x \in X$, there exists an open set $V(x, U)$ containing $x$

\hspace{2mm} such that, $y \in V(x, U)$ implies $x \in W(m, y) \subset U $ for some $m \in \omega$.

\hspace{-8mm}If we strengthen the condition $(G)$ by not allowing the natural number $m$ to vary with $y$, then we say that $X$ satisfies condition $(A)$, that is, for each open set $U$ and $x \in U$, there exists an open set $V(x, U)$ containing $x$ and a natural number $m = m(x, U)$ such that
$x \in W(m, y) \subset U $ for all $y \in V(x, U)$. If each $W(n, x)$ is open, we say that $X$ satisfies open $(G)$ or open $(A)$ respectively.
If $W(n+1, x) \subset W(n, x)$ for each $n \in \omega$, we say that $X$ satisfies decreasing $(G)$ or decreasing $(A)$.
The Collins-Roscoe mechanism has been extensively studied, and a lot of significant results have been obtained. For example,

\hspace{-10mm} \textbf{Theorem 1.1}([1],[2]) The followings are equivalent for a space $X$:

(1) $X$ is metrizable,

(2) $X$ satisfies decreasing open $(A)$,

(3) $X$ satisfies decreasing open $(G)$,

(4) $X$ satisfies decreasing neighborhood $(A)$.

\hspace{-10mm} \textbf{Theorem 1.2}([2],[3]) The followings are equivalent for a space $X$:

(1) $X$ is stratifiable,

(2) $X$ satisfies decreasing $(G)$ and has countable pseudo-character,

(3) $X$ satisfies decreasing $(A)$ and has countable pseudo-character. \\

In [4], E. Van. Douwen first introduced the concept of $D$-space and proved that the finite product of Sorgenfrey lines
is a $D$-space.

\vspace{2mm}
\hspace{-10mm} \textbf{Definition 1.5} A neighborhood assignment for a topological space $(X, \mathcal{T})$
is a function $\phi : X \rightarrow \mathcal{T}$ such that $x \in \phi (x)$. A space $X$ is a $D$-space if
for each neighborhood assignment $\phi$, there is a closed discrete subset $D_{\phi}$ of $X$ such that
$\{\phi (d): d \in D_{\phi}\}$ covers $X$.

 A lot of interesting work on $D$-spaces have been done, in particularly the connections between the generalized metric spaces and $D$-spaces. Borges and Wehrly proved that semi-stratifiable spaces are $D$-spaces [5], Buzyakova showed that every strong $\Sigma$-space is a $D$-space [6], and thus all $\sigma$-spaces are $D$-spaces [7], Arhangel'skii and Buzyakova obtained the interesting result that every space with a point countable base is a $D$-space [8], and so on. For more detail about the work of $D$-spaces, the survey paper [9] written by Gruenhage  is recommended. \\

In [10], Gruenhage developed a new technique based on the earlier work of Fleissner and Stanley, and proved that any space satisfying open $(G)$ is a $D$-space. This established certain connection between the Collins-Roscoe mechanism and $D$-spaces.
Recently, Soukup and Xu proved that the well-ordered ($\alpha$ A), linearly semi-stratifiable space and elastic space are all $D$-spaces. Since these spaces are all well-ordered $(F)$ spaces, they asked that whether well-ordered $(F)$ spaces are $D$-spaces[11]. Considering that well-ordered $(F)$ spaces are monotonically normal paracompact spaces, so this question is a weak version of an old one asked by Borges and Wehrly in [5]: Whether monotonically normal paracompact spaces are $D$-spaces? In the present paper, we prove that well-ordered $(F)$ spaces are $D$-spaces, and thus get more insight in the relationship between the Collins-Roscoe mechanism and $D$-spaces..

Throughout this paper, all spaces are $T_{1}$, and $\omega$ is the first countable ordinal. For other undefined terms we refer the reader to [2], [7] and [9]. \\

\begin{center}{\bf \Large 2.  Well-ordered $(F)$ spaces are $D$-spaces}\end{center}

\par\vspace{6mm}

Recall that, a space $X$ has a family $\mathcal{W}$ satisfying condition $(F)$ if $\mathcal{W} = \{\mathcal{W}(x): x \in X \}$
where each $\mathcal{W}(x)$ consists of subsets of $X$ containing $x$ and

\hspace{4mm} if $x$ belongs to open $U$, then there exists open $V = V(x, U)$ containing $x$

\hspace{4mm} such that $x \in W \subset U $ for some $W \in \mathcal{W}(y)$ whenever $y \in V$.

\hspace{-10mm} We say $X$ satisfies $(F)$ if $X$ has $\mathcal{W}$ satisfying $(F)$. If, in addition, each $\mathcal{W}(x)$ is a chain (well-ordered) under reverse inclusion, then we say that $X$ satisfies chain (well-ordered) $(F)$. Further more, we say $X$ satisfies neighborhood $(F)$ if $X$ has $\mathcal{W}$ satisfying $(F)$, and each element of $\mathcal{W}(x)$ is a neighborhood of $x$.\\

In [12], Stares give a characterization of the decreasing $(G)$ as a strongly monotone normality condition from which we can decide which one of $"x \in V"$ and $"y \in U"$ holds when $H(x, U) \cap H(y, V) \neq \emptyset$. For the sake of completeness, we state it in the following.

\vspace{2mm}
\hspace{-10mm} \textbf{Theorem 2.1}[12] A space satisfies decreasing $(G)$ iff to each point $x \in X$ and open set $U$ containing $x$ we can assign an open set $H(x, U)$ containing $x$ and, for each point $a \in H(x, U)$ we can assign a natural number $n(a, x, U)$ such that if $a \in H(x, U) \cap H(y, V)$ and $n(a, x, U) \leq n(a, y, V)$ then $y \in U$. \\

For the well-ordered $(F)$ space, we have a similar characterization as the decreasing $(G)$ spaces in above theorem. Essentially, the idea of its proof comes from Stares' theorem.

\vspace{2mm}
\hspace{-10mm} \textbf{Theorem 2.2} A space satisfies well-ordered $(F)$ iff to each point $x \in X$ and open set $U$ containing $x$ we can assign an open set $H(x, U)$ containing $x$ and, for each point $a \in H(x, U)$ we can assign an ordinal number $n(a, x, U)$ such that if $a \in H(x, U) \cap H(y, V)$ and $n(a, x, U) \leq n(a, y, V)$ then $y \in U$.

\hspace{-10mm} \textbf{Proof} If $X$ satisfies well-ordered $(F)$, we define $H(x, U) = V(x, U)$ for each $x \in X$ and open set $U$ containing $x$. Let $a \in H(x, U)$, then there is some ordinal number $\alpha$ such that $x \in W(\alpha, a) \subset U$ according to the condition $(F)$. So, we can assign $n(a, x, U) = \alpha$. It is easy to check that $a \in H(x, U) \cap H(y, V)$ and $n(a, x, U) \leq n(a, y, V)$ implies $y \in U$ by the well-ordered property of $\mathcal{W}(a)$.

Conversely, for each $a \in X$ we define $W(\alpha, a) = \{a\} \cap \{y: a \in H(y, V)$ and $n(a, y, V) \geq \alpha$ for some open $V \}$, $\mathcal{W}(a) = \{W(\alpha, a): \alpha \leq \alpha_{a}\}$ where $\alpha_{a} = sup \{n(a, y, V): a \in H(y, V)$ for some open $V \}$, and $\mathcal{W} = \{\mathcal{W}(a): a \in X\}$. Obviously, we have $a \in W(\alpha, a)$ for each $\alpha \leq \alpha_{a}$, and $W(\beta, a) \subset W(\gamma, a)$ whenever $\gamma < \beta \leq \alpha_{a}$.

If $a \in H(x, U)$, we consider the element $W(n(a, x, U), a)$ of $\mathcal{W}(a)$. By the definition of $W(n(a, x, U), a)$, it is obvious that $x \in W(n(a, x, U), a)$. Further more, we shall prove that $W(n(a, x, U), a) \subset U$. For any $y \in W(n(a, x, U), a)$ there are two cases: $y = a$, or $y \neq a$. If $y = a$, then we have $y = a \in H(y, U) \subset U$. Otherwise, $a \in H(y, V)$ and $n(a, y, V) \geq n(a, x, U)$ for some open $V$ according to the definition of $W(n(a, x, U), a)$. This implies that $y \in U$. \\

 The above theorem make us to have the ability to decide which one of "$x \in V$" and $"y \in U"$ holds when $ H(x, U) \cap H(y, V) \neq \emptyset$ in the well-ordered $(F)$ spaces. Next, we will use it to prove our main theorem in this paper.

\vspace{2mm}
\hspace{-10mm} \textbf{Theorem 2.3} Well-ordered $(F)$ spaces are $D$-spaces.

\hspace{-10mm} \textbf{Proof} Let $X$ be a well-ordered $(F)$ space, and $\phi: x \rightarrow \phi (x)$ a neighborhood assignment on $X$. By theorem 2.2, for each $x \in X$ there is an open subset $H(x, \phi (x))$ containing $x$, and for each $a \in H(x, \phi (x))$ there is an ordinal number $n(a, x, \phi (x))$ such that: $a \in H(x, \phi (x)) \cap H(y, \phi (y))$ and $n(a, x, \phi (x)) \leq n(a, y, \phi (y))$ implies $y \in \phi (x)$.

Take an element $x_{0} \in X$, then $\phi(x_{0})\subset X $. If there is some $y \in X \setminus \phi(x_{0})$ such that $x_{0} \in H(y, \phi (y))$, then there is an ordinal number $n(x_{0}, y, \phi (y))$ assigned.
Put $C(x_{0}) = \{y \in X \setminus \phi(x_{0}): x_{0} \in H(y, \phi(y))\}$, and let $m(x_{0}) = min \{n(x_{0}, y, \phi (y)): y \in C(x_{0})\}$, then there is some $y_{0} \in C(x_{0})$ such that $m(x_{0}) = n(x_{0}, y_{0}, \phi (y_{0}))$. For each $y \in C(x_{0})$, we have $n(x_{0}, y, \phi (y)) \geq m(x_{0}) = n(x_{0}, y_{0}, \phi (y_{0}))$. Notice that $x_{0} \in H(y, \phi (y)) \cap H(y_{0}, \phi (y_{0}))$, by theorem 2.2 we conclude that $y \in \phi (y_{0})$. Therefore, $C(x_{0}) \subset \phi(y_{0})$ holds. Let $x_{1} = y_{0}$, then $C(x_{0}) \subset \phi(x_{1})$ and for each $y \in X \setminus (\phi(x_{0})\cup \phi(x_{1}))$, we have $x_{0} \notin H(y, \phi(y))$.

If for each $y \in X \setminus \phi(x_{0})$ we have $x_{0} \notin H(y, \phi(y))$, take an element of $X \setminus \phi(x_{0})$ as $x_{1}$. Then, $x_{0} \notin H(y, \phi(y))$ also holds for each $y \in X \setminus (\phi(x_{0})\cup \phi(x_{1}))$.

Suppose that $\alpha \geq 1$ is an ordinal, and for each $\beta < \alpha$ we have selected the element $x_{\beta} \in X$ such that
$y \in X \setminus (\cup \{\phi(x_{\delta}): \delta \leq \beta \})$ implies $x_{\gamma} \notin H(y, \phi(y))$ for each $\gamma < \beta$.
Next, if $X \setminus (\cup \{\phi(x_{\delta}): \delta < \alpha \}) \neq \emptyset$, we shall take some element of $X \setminus (\cup \{\phi(x_{\delta}): \delta < \alpha \})$ as $x_{\alpha}$ such that

\hspace{6mm} $y \in X \setminus (\cup \{\phi(x_{\delta}): \delta \leq \alpha \})$ implies $x_{\gamma} \notin H(y, \phi(y))$ for each $\gamma < \alpha$. \ \ \ \ \ \ $(\ast)$

First, if $\alpha$ is a limit ordinal, we take an element of $X \setminus (\cup \{\phi(x_{\delta}): \delta < \alpha \})$ as $x_{\alpha}$. For each $\gamma < \alpha$, we have $\gamma < \gamma + 2 < \alpha$. If $y \in X \setminus (\cup \{\phi(x_{\delta}): \delta \leq \alpha \})$, then $y \in X \setminus (\cup \{\phi(x_{\delta}): \delta \leq \gamma + 2 \})$. By the inductive assumption, we claim that $x_{\gamma} \notin H(y, \phi (y))$.

Second, if $\alpha$ is a success ordinal, and $y \in X \setminus (\cup \{\phi(x_{\delta}): \delta \leq \alpha - 1\})$. By the inductive assumption, we have that

\hspace{24mm} $x_{\gamma} \notin H(y, \phi (y))$ for each $\gamma < \alpha - 1$. \ \ \ \ \ \ \  \ \ \ \ \ \ \   $(\ast \ast)$

  If $y \in X \setminus (\cup \{\phi(x_{\delta}): \delta \leq \alpha - 1 \})$ implies $x_{\alpha - 1} \notin H(y, \phi (y))$, we take an element of $X \setminus (\cup \{\phi(x_{\delta}): \delta \leq \alpha - 1 \})$ as $x_{\alpha}$. Obviously, for each $y \in X \setminus (\cup \{\phi(x_{\delta}): \delta \leq \alpha \})$, we have $x_{\alpha - 1} \notin H(y, \phi (y))$. With the fact $(\ast \ast)$, we conclude that $x_{\alpha}$ satisfies the condition $(\ast)$.

Otherwise, there is some $y \in X \setminus (\cup \{\phi(x_{\delta}): \delta \leq \alpha - 1 \})$ such that $x_{\alpha - 1} \in H(y, \phi (y))$ for which the corresponding ordinal number is $n(x_{\alpha - 1}, y, \phi (y))$. Denote $C(x_{\alpha - 1}) = \{y \in X \setminus (\cup \{\phi(x_{\delta}): \delta \leq \alpha -1 \}): x_{\alpha - 1} \in H(y, \phi (y))$, and let $m(x_{\alpha - 1}) = min\{n(x_{\alpha - 1}, y, \phi (y)): y \in C(x_{\alpha - 1}) \}$, then there is some $y_{\alpha} \in C(x_{\alpha - 1})$ such that $x_{\alpha - 1} \in H(y_{\alpha}, \phi (y_{\alpha}))$ and $n(x_{\alpha - 1}, y_{\alpha}, \phi (y_{\alpha})) = m(x_{\alpha - 1})$. If $y \in C(x_{\alpha - 1})$, then $x_{\alpha - 1} \in H(y, \phi (y))$ and $m(x_{\alpha - 1}) \leq n(x_{\alpha - 1}, y, \phi (y))$, and thus $n(x_{\alpha - 1}, y_{\alpha}, \phi (y_{\alpha})) \leq n(x_{\alpha - 1}, y, \phi (y))$. It follows that $y \in \phi (y_{\alpha})$. Hence, we have that $C(x_{\alpha - 1}) \subset \phi (y_{\alpha})$. Let $x_{\alpha}= y_{\alpha}$, then for $y \in X \setminus (\cup \{\phi(x_{\delta}): \delta \leq \alpha \})$ implies that $x_{\alpha - 1} \notin H(y, \phi (y))$. Notice the fact $(\ast \ast)$, we claim that
$x_{\alpha}$ satisfies the condition $(\ast)$.

By the transfinite induction, we can get a subset $D_{\phi} = \{x_{\alpha}: \alpha < \kappa\}$ of $X$ such that each $x_{\alpha}$ satisfies the condition $(\ast)$ and $\cup \{ \phi(x_{\alpha}): \alpha < \kappa \} = X$, where $\Gamma$ is an indexed set of ordinals.

Let $x \in X$, and $\alpha$ be the smallest ordinal such that $x \in \phi (x_{\alpha})$. If $\alpha$ is a limit ordinal, we claim that $x_{\beta} \notin H(x, \phi(x))$ for each $\beta < \alpha$. In fact, suppose $x_{\beta_{0}} \in H(x, \phi(x))$ for some $\beta_{0} < \alpha$, then we have $x \in X \setminus (\cup \{\phi (x_{\delta}): \delta \leq \beta_{0}\})$. By the selection of $x_{\beta_{0}+1}$, we know that $x \in C(x_{\beta_{0}}) \subset \phi(x_{\beta_{0} + 1})$, a contradiction with $\beta_{0} + 1  < \alpha$. If $\alpha$ is a successor ordinal, by condition  $(\ast \ast)$ we have that $x_{\beta} \notin H(x, \phi (x))$ for each $\beta < \alpha - 1$.

For each $\beta > \alpha$, we know that $x_{\beta} \in X \setminus (\cup \{\phi(x_{\delta}): \delta \leq \alpha \})$, and thus $x_{\beta} \notin \phi (x_{\alpha})$. Therefore, we find an open neighborhood $H(x, \phi(x)) \cap \phi (x_{\alpha})$ of $x$ such that $|(H(x, \phi(x)) \cap \phi (x_{\alpha})) \cap D_{\phi}| \leq 2$. Since $X$ is a $T_{1}$-space, we conclude that $D_{\phi}$ is a closed discrete subset of $X$, and thus $X$ is a $D$-space. \\

The converse of the above theorem is not true. There is an example of $D$-space which fails to have well-ordered $(F)$.

\vspace{2mm} \hspace{-10mm} \textbf{Example 2.4} Let $X = R \cup (\cup \{Q \times \{ \frac{1}{n}\}: n \in N \})$.
The topology is defined as following: The point of $X \setminus R$ is isolated.
For $x \in R$, the element of its neighborhood base is $\{ x \}  \cup (\cup \{([a_{x,n},x) \cap Q) \times \{ \frac{1}{n} \}: n \geq m \})$.

Lin proved that $X$ is k-semi-stratifiable but not normal[13]. Since the k-semi-stratifiable space is obviously semi-stratifiable,
$X$ is a $D$-space[11]. On the other hand, both decreasing $(G)$ space and well-ordered $(F)$ space are monotonically normal[12, 14], so $X$ does not satisfy decreasing $(G)$ and well-ordered $(F)$. \\

\begin{center}{\bf \Large 3. Further discussion}\end{center}

\par\vspace{6mm}

In [2], Collins et. provided us an example of a monotonically normal space having $\mathcal{W}$ satisfying
chain $(F)$ but not metacompact, that is, the ordinal space $[0, \omega_{1})$. Since monotonically normal $D$-space is paracompact[5], we conclude that $[0, \omega_{1})$ is not a $D$-space. This fact tells us that the chain $(F)$ spaces need not to be $D$-spaces, although the well-ordered $(F)$ spaces are $D$-spaces.

When Borges and Wehrly proved that  monotonically normal $D$-space is paracompact in [5], they asked the question "Whether every monotonically normal paracompact space is $D$-space?" Until now, it is still open. In [2], Collins et. indicated that both well-ordered (F) space and chain neighborhood (F) space are paracompact (Moody et. gave them an interesting name "the unified paracompactness theorem" [15]). In fact, they are both monotonically normal paracompact spaces, so the following question is of interesting.

\vspace{2mm} \hspace{-10mm} \textbf{Question 3.1} Whether every neighborhood chain $(F)$ space is a $D$-space? \\

In [16], the well-known Smirnov metrization theorem is shown, i.e., a space is metrizable iff it is paracompact and locally metrizable. Gao generalized this result to the decreasing $(G)$ space, and proved that a space satisfies decreasing $(G)$ iff it is a paracompact and local satisfies decreasing $(G)$. Similar as the proof of theorem 2.6 in [17], we can prove that this situation is also true for the well-ordered $(F)$ spaces.

\vspace{2mm} \hspace{-10mm} \textbf{Theorem 3.2} A space satisfies well-ordered $(F)$ iff it is a paracompact and local satisfies well-ordered $(F)$.

In [18], Moody and Roscoe proved that a monotonically normal space is acyclic monotonically normal iff it can be covered by a collection of open acyclic monotonically normal subspaces. That is to say a monotonically normal space is chain $(F)$ iff it can be covered by a collection of open chain $(F)$ subspaces. For the well-ordered $(F)$ space, it is reasonable to ask the following question.

\vspace{2mm} \hspace{-10mm} \textbf{Question 3.3} Whether a monotonically normal space is well-ordered $(F)$ iff it can be covered by a collection of open well-ordered $(F)$ subspaces.


\hspace{2cm}
\begin{center}{\large \bf  REFERENCES}\end{center}
\hspace{0.5cm}


\end{document}